\documentclass[a4paper,12pt]{amsart}

\newtheorem{theo}{Theorem}[section]

\newtheorem{prop}[theo]{Proposition}
\newtheorem{lem}[theo]{Lemma}

\newtheorem{cor}[theo]{Corollary}

\newtheorem{rema}[theo]{Remark}
\newtheorem{remas}[theo]{Remarks}

\def \Romannumeral #1 {\expandafter\uppercase\expandafter {\romannumeral #1} }
\def \br {{\rm{Br\,}}}

\def \P {{\bf P}}

\def \Ga {{\Gamma}}

\def \pic {{\rm {Pic\,}}}

\def \dim {{\rm{dim\,}}}
\def \ns {{\rm{NS\,}}}
\def \Hom {{\rm {Hom}}}

\def \Z {{\bf Z}}
\def \Q {{\bf Q}}

\def \C {{\bf C}}

\def \G {{\bf G}_m}
\def \Ga {{\bf G}_a}

\def \Het {H_{\mbox{\scriptsize\'et}}}

\def\smallsquare{\vbox{\hrule\hbox{\vrule height 1 ex\kern 1 ex\vrule}\hrule}}
\def\enddem{\hfill \smallsquare\vskip 3mm}

\usepackage{amscd}
\usepackage{amssymb}
\usepackage{amsmath}

\everymath{\displaystyle}

\DeclareFontFamily{U}{wncy}{}
\DeclareFontShape{U}{wncy}{m}{n}{%
   <5>wncyr5%
   <6>wncyr6%
   <7>wncyr7%
   <8>wncyr8%
   <9>wncyr9%
   <10>wncyr10%
   <11>wncyr10%
   <12>wncyr6%
   <14>wncyr7%
   <17>wncyr8%
   <20>wncyr10%
   <25>wncyr10}{}
\DeclareMathAlphabet{\cyrille}{U}{wncy}{m}{n}

\title{Prime-to-$p$ \'etale covers of algebraic groups and homogeneous spaces}
\author{Michel Brion and Tam\'as Szamuely}
\address{Universit\'e de Grenoble I, D\'epartement de Math\'ematiques, Institut Fourier, UMR 5582 du CNRS,
38402 Saint-Martin d'H\`eres Cedex, France}
\email{Michel.Brion@ujf-grenoble.fr}
\address{Alfr\'ed R\'enyi Institute of Mathematics, Hungarian Academy of Sciences, Re\'altanoda utca 13-15, H-1053 Budapest, Hungary}
\email{szamuely.tamas@renyi.mta.hu}

\date{\today}
\begin{document}
\maketitle

\markright{PRIME-TO-$p$ \'ETALE COVERS}

\section{Introduction}

By a classical result of Schreier \cite{schreier}, the fundamental group of a connected and locally connected topological group $G$ is commutative. If moreover $G$ is (semi-)locally simply connected, then every Galois cover $\phi:\, Y\to G$ carries a group structure for which $\phi$ is a homomorphism. Thus $G$ is the quotient of $Y$ by an abelian normal subgroup.

The first part of the following proposition states an analogue of this result in algebraic geometry. The second part gives a bound on the number of topological generators of the prime-to-$p$ fundamental group. To state it, we need to introduce some notation. Recall that by Chevalley's theorem $G$ is an extension of an abelian variety $A$ by a linear algebraic group $G_{\rm aff}$ (see \cite{bsu}, \cite{chev}, \cite{conrad}, \cite{ros}). Denote by $g$ the dimension of $A$ and by $r$ the rank of $G_{\rm aff}$ (which is by definition the dimension of a maximal torus). Furthermore, denote by $\Z_{(p')}$ the direct product of the rings $\Z_\ell$ for $\ell \neq p$.

\begin{prop}\label{pi1p} Let $G$ be a connected algebraic group over an algebraically closed field of characteristic $p\geq 0$.
\begin{enumerate}
\item[a)] Every \'etale Galois cover $Y\to G$ of degree prime to $p$ carries the structure of a central isogeny. In particular, the maximal prime-to-$p$ quotient $\pi_1^{(p')}(G)$ of the \'etale fundamental group of $G$ is commutative.
\item[b)] The group $\pi_1^{(p')}(G)$ is a quotient of $\Z_{(p')}^{2 g + r}$.
\end{enumerate}
\end{prop}

By convention, all Galois covers are assumed to be connected.
Statement $a)$ was known in the special cases where $G$ is an abelian variety (Lang--Serre \cite{langserre}) or $G$ is affine (Miyanishi \cite{miyanishi}, Magid \cite{magidcov}).
Note that restricting to prime-to-$p$ covers is crucial because if one allows $p$-covers in characteristic $p>0$, the statement fails disastrously already for the additive group $\Ga$ (Raynaud \cite{raynaud}).

In part $b)$ the bound is sharp (take the direct product of an abelian variety and a torus). In the case when $G$ is commutative it follows from Theorem \ref{pi1p} and (\cite{proalg}, \S 6.4, Cor. 4) that $\pi_1^{(p')}(G)$ is moreover a free $\Z_{(p')}$-module. However, in the non-commutative case it may contain torsion. For instance, $\pi_1({\rm SO}(n))$ in characteristic 0 is of order 2. (Of course, in this case the bound on the number of generators given by the proposition is very far from the truth.)

We also prove a generalization to homogeneous spaces. Bounding the topological rank is more difficult than in the group case, so we state this result as a theorem.

\begin{theo}\label{eff} Keep the assumptions of the proposition above, and let $X$ be a variety on which $G$ acts transitively with connected stabilizers.
\begin{enumerate}
\item[a)]  Every \'etale Galois cover of $X$ of degree prime to $p$ is of the form $\tilde{G}/\tilde{H}$ for a central isogeny
$\tilde{G} \to G$ and a lift $\tilde{H}$ of the stabilizer $H$ of a point of $X$ in $\tilde{G}$.
\item[b)] The group $\pi_1^{(p')}(X)$ is a quotient of $\Z_{(p')}^{2(g-g_H)+ r- r_H}$,
where $g_H$ and $r_H$ are the dimension of the abelian variety quotient and the rank of the subgroup $H$, respectively.
\end{enumerate}
\end{theo}
In work in progress, M. Borovoi and C. Demarche \cite{bd} have obtained an extension of statement $a)$ to non-connected $H$. In this case $\tilde H$ is a finite index subgroup of the preimage of $H$ in $\tilde G$.
However, statement $b)$ fails when the isotropy group $H$ is not connected.
For example, if $\Gamma$ is a finite group of order prime to $p$ and $\Gamma \hookrightarrow {\rm GL}(n)$ a faithful representation, then $\pi_1^{(p')}({\rm GL}(n)/\Gamma)$ is just $\Gamma$.

We also present an application of Proposition \ref{pi1p} which was our original motivation for studying the problem. It concerns commutative algebraic groups.

\begin{theo}\label{berkex} Let $G$ be a connected commutative algebraic group over an algebraically closed field $k$. For $n$ prime to the characteristic there is an exact sequence
$$
0\to \pic G/n\/\pic G \to \Hom(\Lambda^2 ({}_nG), \mu_n)\to {}_n\br G\to 0
$$
where $\pic G$ and $\br G$ denote the Picard and Brauer groups of $G$, respectively, and ${}_nG$ stands for the $n$-torsion part of $G$.
\end{theo}

The case where $G$ is an abelian variety was established by Berkovich \cite{berk}. The proof in the general case uses Proposition \ref{pi1p}  as an input, in the identification of the group $\Hom({}_nG, \mu_n)$
with the \'etale cohomology group $H^1(G,\mu_n)$. However, this identification can also be derived using a more elementary method; see the conclusion of Section \ref{secberkex}.

The following almost immediate corollary of Theorem \ref{berkex} makes a classical formula
of Grothendieck (\cite{gb}, corollaire 3.4) explicit in the case of a connected commutative
algebraic group $G$.

\begin{cor}\label{corberkex} Let $\ell$ be a prime number different from the characteristic of $k$. The $\ell$-primary torsion part of $\br G$ is isomorphic to $(\Q_\ell/\Z_\ell)^{(2g+r)(2g+r-1)/2 - \rho}$, where $g$ and $r$ are defined as above, and $\rho$ is the rank of the N\'eron--Severi group of $A$, the abelian variety quotient of $G$.
\end{cor}

Theorem \ref{berkex} and Corollary \ref{corberkex} will be proven in Section \ref{secberkex}.

\medskip

We thank Mikhail Borovoi, Antoine Chambert-Loir, Cyril Demarche, H\'el\`ene Esnault, Ofer Gabber and Jakob Stix for very helpful remarks. The second author was partially supported by OTKA grant No. NK81203.

\section{Algebraic groups}

In this section we prove Proposition \ref{pi1p}. Our first proof of statement $a)$ was quite involved, using, among other ingredients, a fairly recent result of Orgogozo \cite{org} on compatibility of the prime-to-$p$ fundamental groups with products. As H\'el\`ene Esnault pointed out, this result together with a modification of the classical argument of Lang and Serre for abelian varieties suffices for the proof of the statement. We briefly recall the argument (see the original paper \cite{langserre} or \cite{fg}, Theorem 5.6.10 and Remark 5.6.11 for more details).

Denote by $m_G:\, G\times G\to G$ the group operation of $G$. Given an \'etale Galois cover $Y\to G$ of degree prime to $p$,  consider the fibre
product square
$$
\begin{CD}
Y'=(G\times G)\times_GY @>>> Y \\
@VVV @VV{\phi}V \\
G\times G @>{m_G}>> G.
\end{CD}
$$
One checks that $Y'$ is also connected, hence $Y'\to G\times G$ is a Galois \'etale cover with group $\Gamma:={\rm Aut}(Y|G)$. Since the prime-to-$p$ fundamental group is compatible with direct products (Orgogozo \cite{org}, Cor. 4.9),  there exist \'etale Galois covers $Z_1\to G$ and $Z_2\to G$, with groups $\Gamma_1$
and $\Gamma_2$, respectively, such that $Y'$ is a quotient of the direct product $Z_1\times
Z_2\to G\times G$. In particular, there is a normal subgroup $\Delta\subset \Gamma_1\times \Gamma_2$ with
$\Gamma\cong (\Gamma_1\times \Gamma_2)/\Delta$. Replacing the $Z_i$ by their quotients by the actions of the
$\Gamma_i\cap \Delta$ we may assume $\Delta\cap \Gamma_i=\{1\}$ for $i=1,2$. Thus the $\Gamma_i$ may be identified
with normal subgroups of $\Gamma$ that generate $\Gamma$ and centralize each other. But by
construction the base change $Y_2'\to \{1\}\times G$ of $Y'$ by the map $\{1\}\times G\to G\times G$ must be an \'etale Galois cover with group $\Gamma_2$, and
similarly for the base change $Y_1'\to G\times\{1\}$. It follows that $\Gamma=\Gamma_1=\Gamma_2$, and $\Gamma$ is commutative. Moreover, we have isomorphisms $Z_i\cong Y_i'\cong Y$ for $i=1,2$. This
yields a map $Y\times Y\to Y'$, whence also a map $m_Y:\, Y\times Y\to Y$ by composing
with the projection $Y'\to Y$. Fix a point $1_Y$ of $Y$ in the fibre above the neutral element of $G$. Modifying
$m_Y$ by an automorphism of $Y$ if necessary we may assume $m_Y(1_Y, 1_Y)=1_Y$. A standard argument (\cite{fg}, pp. 180--181) then shows that $m_Y$ equips $Y$ with the structure of an algebraic group.\medskip

We now turn to the proof of part $b)$ and take up the notation introduced before the statement of the proposition.
We already know  that $\pi_1^{(p')}(G)$ is commutative.
In the case when $G$ is an abelian variety of dimension $g$ it is moreover free over
$\Z_{(p')}$ of rank $2g$ (see e.g. \cite{fg}, Theorem 5.6.10 and Remark 5.6.11).
Hence it will be enough to prove the following: the natural map
$\alpha_*:\,\pi_1^{(p')}(G)\to \pi_1^{(p')}(A)$ induced by the projection
$\alpha: \, G \to A$ is surjective and its kernel is a quotient of $\Z_{(p')}^r$.

Denote by $\Q_{(p')}/\Z_{(p')}$ the direct sum of the torsion abelian groups
$\Q_\ell/\Z_\ell$ for $\ell\neq p$. Since $\pi_1^{(p')}(G)$ is commutative,
there is a perfect pairing
$$
\pi_1^{(p')}(G)\times \Het^1(G, \Q_{(p')}/\Z_{(p')})\to \Q/\Z
$$
obtained by passing to the limit from the perfect pairings
$$
\pi_1^{(p')}(G)/m \pi_1^{(p')}(G)\times \Het^1(G, \Z/m\Z)\to \Z/m\Z
$$
for $m$ prime to $p$.

The Leray spectral sequence
$$
\Het^{p}(A, {\bf R}^q\alpha_*\Q_{(p')}/\Z_{(p')})\Rightarrow \Het^{p+q}(G, \Q_{(p')}/\Z_{(p')})
$$
associated with the projection $\alpha$ yields an exact sequence
\begin{align*}
0\to\Het^{1}(A, \alpha_*\Q_{(p')}/\Z_{(p')}) &\to \Het^{1}(G, \Q_{(p')}/\Z_{(p')})\to
\\ &\to \Het^{0}(A, {\bf R}^1\alpha_*\Q_{(p')}/\Z_{(p')}).
\end{align*}

\noindent As $\alpha$ has connected fibres,
$\Het^{1}(A, \alpha_*\Q_{(p')}/\Z_{(p')})\cong \Het^{1}(A, \Q_{(p')}/\Z_{(p')})$,
which shows that the map between the first two terms in the above exact sequence
is the $\Q/\Z$-dual of the map $\alpha_*:\,\pi_1^{(p')}(G)\to \pi_1^{(p')}(A)$.
It follows that the latter map is surjective.

The stalk of the \'etale sheaf ${\bf R}^1\alpha_*\Q_{(p')}/\Z_{(p')}$
at a geometric point can be identified with $\Het^{1}(G_{\rm aff}, \Q_{(p')}/\Z_{(p')})$.
Pick a maximal torus $T\subset G_{\rm aff}$. Miyanishi proves in (\cite{miyanishi}, Lemma 3)
that the natural map $\pi_1^{(p')}(T)\to \pi_1^{(p')}(G_{\rm aff})$ is surjective
(see also Lemma \ref{maxtorus} below for a more general statement). Since
$\pi_1^{(p')}(T)\cong \pi_1^{(p')}(\G)^r\cong \Z_{(p')}^r$,
this bounds the number of generators of $\pi_1^{(p')}(G_{\rm aff})$ by $r$.
But since $\alpha$ is smooth, the sheaf
${\bf R}^1\alpha_*\Q_{(p')}/\Z_{(p')}$ is locally constant and thus
$\Het^{0}(A, {\bf R}^1\alpha_*\Q_{(p')}/\Z_{(p')})$ becomes a subgroup of
$\Het^{1}(T, \Q_{(p')}/\Z_{(p')})\cong (\Q_{(p')}/\Z_{(p')})^r$.
By taking $\Q/\Z$-duals we obtain that the kernel of
$\alpha_*:\,\pi_1^{(p')}(G)\to \pi_1^{(p')}(A)$ is indeed a quotient of $\Z_{(p')}^r$.

\begin{remas}\rm ${}$

\noindent 1. In the course of the above proof we have in fact constructed a homotopy exact sequence
$$
\pi_1^{(p')}(G_{\rm aff})\to\pi_1^{(p')}(G)\to \pi_1^{(p')}(A)\to 0.
$$
 The paper \cite{bd} contains a generalization of this exact sequence to homogeneous spaces which is constructed using our Theorem \ref{eff} $a)$. Note that a general homotopy exact sequence of prime-to-$p$ fundamental groups for morphisms with non-proper fibres is not known at present.

\noindent 2. Jakob Stix has shown us the following very quick argument for the commutativity of $\pi_1^{(p')}(G)$.  The point is that by a well-known reasoning $\pi_1^{(p')}(G)$ is a group object in the category of groups. In more detail: by compatibility of the functor $\pi_1^{(p')}$ with products (\cite{org}, Cor. 4.9 again) the multiplication map of $G$ induces a map $m: \pi_1^{(p')}(G)\,\times\, \pi_1^{(p')}(G)\to \pi_1^{(p')}(G)$ that satisfies the group axioms by functoriality, and moreover $m(1,g)=m(g,1)=g$ for all $g\in \pi_1^{(p')}(G)$. The calculation
$$
hg = m(1,h)m(g,1) = m((1,h)(g,1)) = m(g,h) $$ $$= m((g,1)(1,h)) = m(g,1)m(1,h) = gh
$$
shows that this group law is commutative and coincides with the usual group law of $\pi_1^{(p')}(G)$.
\end{remas}

\section{Homogeneous spaces}
\label{sechom}

In this section we prove Theorem \ref{eff}.
We begin the proof of part $a)$ with some standard observations and reductions. Let $X$ be a variety on which a connected algebraic group $G$ acts transitively.
Choose a point $x \in X$ and denote by $H \subset G$ its isotropy subgroup scheme. Then
$X \cong G/H$ and the natural map $G/H_{\rm red} \to G/H$ is finite, surjective and purely inseparable,
where $H_{\rm red}$ denotes the reduced part of $H$. It follows that the induced map
$\pi_1(G/H_{\rm red}) \to \pi_1(G/H)$ is an isomorphism. So we may assume throughout that
$H$ is is an algebraic group; moreover, as in the statement of the theorem, we assume it is connected.

Assume now $f : \tilde{X} \to X$ is an \'etale Galois cover with group $\Gamma$ of order prime to $p$. Consider the orbit map $\varphi : G \to X$, $g \mapsto g \cdot x$, and form the cartesian square
$$
\CD
\tilde{G} @>{\tilde{\varphi}}>> \tilde{X} \\
@V{\pi}VV @V{f}VV \\
G @>{\varphi}>> X.\\
\endCD
$$
Since $\varphi$ is a flat (hence open) surjection with connected fibres, so is $\tilde\varphi$, so we infer from the connectedness of $\tilde X$ that $\tilde G$ is connected as well.
By Proposition \ref{pi1p} $a)$ the map $\pi$ is therefore a central isogeny with kernel $\Gamma$. The algebraic subgroup
$\pi^{-1}(H)\subset \tilde{G}$ is a disjoint union of fibres of $\tilde{\varphi}$; denote by $\tilde{H}$ its neutral component. Since $H$ is connected by assumption, the map $\pi|_{\tilde H}:\,\tilde{H}\to H$ must be an isomorphism. Moreover, the map $\tilde{\varphi}$ is invariant under $\tilde{H}$ since so is $f \tilde{\varphi} = \varphi \pi\/$ and $f$ is finite \'etale (see e.g. \cite{fg}, Corollary 5.3.3).
Thus $\tilde\varphi$ factors through a morphism $\tilde{G}/\tilde{H} \to \tilde{X}$, which is an
isomorphism since so is the morphism $G/H \to X$. This completes the proof of statement $a)$.\medskip

We now turn to the rank estimate of part $b)$. Given $X = G/H$ as in the statement, we may assume without harm that $G$ acts faithfully on $X$. Indeed, the (scheme-theoretic) kernel $K$ of the action is contained in $H$. Then $G/K$ acts faithfully on $X$ with connected stabilizer $H/K$.
 By the lemma on p. 154 of \cite{matsumura} (see also \cite{bsu}, Corollary 2.1.9) the quotient $H/K$ is  affine. Thus the dimension $g_H$ of the maximal abelian quotient of $H$ equals that of the reduced neutral component of $K$, and hence the dimension of the maximal abelian quotient of $G/K$ is $g-g_H$.  Thus we may replace $G$ by $G/K$ and $H$ by $H/K$ and have a natural map $\alpha_X : X = G/H \to G/G_{\rm aff} = A$
(which is in fact the Alban\-ese morphism of $X$, but we shall not need this). \smallskip

\noindent {\em Case I: the base field $k$ has characteristic $0$.} The fundamental group of $X$ being topologically finitely generated (\cite{sga7}, Exp. II), it does not change by extensions of algebraically closed fields by the argument of (\cite{fg}, p. 186; see also \cite{org}, Theorem 4.11). Thus we reduce to the case $k = \C$ where the algebraic fundamental group is the profinite completion of the topological fundamental group $\pi_1^{\rm top}(X)$ (\cite{sga1}, Exp. XII, Cor. 5.2). It is then enough to see that $\pi_1^{\rm top}(X)$ is a quotient of $\Z^{2 g + r - r_H}$.
The fibration $\alpha_X$ with fibre $G_{\rm aff}/H$ yields an exact sequence
$$
\pi_1^{\rm top}(G_{\rm aff}/H) \to \pi_1^{\rm top}(X) \to \pi_1^{\rm top}(A) \cong \Z^{2g}.
$$
Thus it remains to show that $\pi_1^{\rm top}(G_{\rm aff}/H)$ is a quotient of $\Z^{r - r_H}$.

Choose a maximal torus $T_H \subset H$ and then a maximal torus $T \subset G_{\rm aff}$ containing
$T_H$. Then we have a commutative diagram
$$
\CD
T @>>> G_{\rm aff} \\
@VVV @VVV \\
T/T_H @>>> G_{\rm aff}/H \\
\endCD
$$
and hence a similar commutative diagram of induced homomorphisms of
topological fundamental groups.

Let $B$ be a Borel subgroup of $G_{\rm aff}$ containing $T$. The flag variety $G_{\rm aff}/B$
is simply connected, and therefore the natural map
$\pi_1^{\rm top}(B) \to \pi_1^{\rm top}(G_{\rm aff})$ is surjective. Also, $B = T\times  U$ as a variety,
where $U \subset B$ denotes the unipotent part.
Since $U$ is isomorphic to an affine space, the natural map
$\pi_1^{\rm top}(T) \to \pi_1^{\rm top}(B)$
is surjective as well, whence the surjectivity of the composite map $\pi_1^{\rm top}(T) \to \pi_1^{\rm top}(G_{\rm aff})$ (note that Miyanishi's argument mentioned in the previous section is an algebraic variant of this reasoning).

As $G_{\rm aff} \to G_{\rm aff}/H$ is a quotient map with connected fibres, the map $\pi_1^{\rm top}(G_{\rm aff}) \to \pi_1^{\rm top}(G_{\rm aff}/H)$ is also surjective, hence
 so is the map $\pi_1^{\rm top}(T/T_H) \to \pi_1^{\rm top}(G_{\rm aff}/H)$.
But $T/T_H$ is a torus of rank $r - r_H$, therefore
$\pi_1^{\rm top}(T/T_H) \cong \Z^{r - r_H}$, which completes the proof of the characteristic 0 case.\medskip

\noindent {\em Case II: $k$ has characteristic $p>0$.} The argument will be a variant of the above but will use more background from the structure theory of  algebraic groups.

Recall from (\cite{dg}, III.3.8 and \cite{ros}, Corollary 5, p. 440) that
$$
G = G_{\rm aff} G_{\rm ant},
$$
where $G_{\rm ant} \subset G$ denotes the largest anti-affine subgroup
(that is, the largest subgroup with $\mathcal{O}(G_{\rm ant}) = k$); moreover, $G_{\rm ant}$ is connected and central in $G$. In positive characteristic $G_{\rm ant}$ is in fact a semi-abelian variety (\cite{aa}, Proposition 2.2).

As before, choose a maximal torus $T_H \subset H$ and then a maximal torus $T \subset G_{\rm aff}$ containing $T_H$. Concerning $TG_{\rm ant}$ we shall prove below the following key proposition.

\begin{prop}\label{tgant} The natural map
$$
\pi_1(T G_{\rm ant})^{(p')} \to \pi_1(G)^{(p')}
$$
is surjective.
\end{prop}

Assuming the proposition, consider the commutative diagram
$$
\CD
T G_{\rm ant} @>>> G \\
@VVV @VVV \\
(T G_{\rm ant})/T_H @>>> G/H. \\
\endCD
$$
The right vertical map induces a surjection on fundamental groups because the pullback of any connected \'etale cover of $G/H$ to $G$ remains connected, as in the proof of statement $a)$. Proposition \ref{tgant} and the induced diagram on prime-to-$p$ fundamental groups then imply the surjectivity of the map
$$
\pi_1((T G_{\rm ant})/T_H)^{(p')} \to \pi_1(G/H)^{(p')}.
$$

Observe now that $TG_{\rm ant}$ is a semi-abelian variety which is an extension of an abelian variety of dimension $g$ by a torus of rank $r$. To prove this fact, remark first that $T$ contains the maximal
torus of the semi-abelian variety $G_{\rm ant}$, since that torus is central
in $G_{\rm aff}$. Thus $T G_{\rm ant}$ is also a semi-abelian variety, with maximal
torus $T$ and abelian quotient
$$
T G_{\rm ant}/T \cong G_{\rm ant}/(G_{\rm ant} \cap T).
$$
But $G_{\rm ant} \cap T$ is an affine subgroup scheme of $G_{\rm ant}$
containing $(G_{\rm ant})_{\rm aff}$ (the maximal torus of $G_{\rm ant}$),
hence the quotient $(G_{\rm ant} \cap T)/(G_{\rm ant})_{\rm aff}$ is finite.
Likewise, the quotient $(G_{\rm ant} \cap G_{\rm aff})/(G_{\rm ant})_{\rm aff}$
is finite. It follows that the abelian quotient of $T G_{\rm ant}$ is
isogenous to $G_{\rm ant}/(G_{\rm aff} \cap G_{\rm ant}) \cong G/G_{\rm aff} = A$.

We conclude that $(T G_{\rm ant})/T_H$ is
an extension of an abelian variety of dimension $g$ by a torus of dimension $r - r_H$,
and hence its prime-to-$p$ fundamental group is a quotient of $\Z_{(p')}^{2 g + r - r_H}$ by the results of the previous section. This finishes the proof of case II modulo Proposition \ref{tgant}.\medskip

The proof of Proposition \ref{tgant} will require a series of lemmas.

\begin{lem}\label{Gators} If $f : X \to Y$ is a torsor under a connected unipotent group,
then the natural map $\pi_1(X) \to \pi_1(Y)$ is surjective.\end{lem}

\begin{dem} By writing the unipotent group as an iterated extension of copies of $\Ga$ we reduce by induction to the case where $f$ is a $\Ga$-torsor.
Then $f$ has a compactification by the projective line bundle
${\bar{f} : X \times^{\Ga} \P^1 \to Y}$, where $\Ga$ acts on $\P^1$ by translations.
The complement of $X$ in $X \times^{\Ga} \P^1$ is a divisor which yields a section of
$\bar{f}$. So the assertion follows from \cite{sga1}, Expos\'e XIII, Proposition 4.1 and
Example 4.4.\end{dem}

Returning to our situation, we choose a Borel subgroup $B$ of $G$ (or equivalently, of $G_{\rm aff}$) containing $T$.
Then $B G_{\rm ant}$ is a connected solvable subgroup of $G$, since $G_{\rm ant}$
is connected and central in $G$.

\begin{lem}\label{maxsolv} The natural map
$$
\pi_1(B G_{\rm ant})^{(p')} \to \pi_1(G)^{(p')}
$$
is surjective.
\end{lem}

\begin{dem} Choose another Borel subgroup $B^-$ such that the intersection
$B \cap B^-$ is of smallest dimension.
Specifically, consider the unipotent radical $R_u(G) = R_u(G_{\rm aff})$.
Then $G_{\rm aff}/R_u(G)$ is a connected reductive group with Borel subgroup
$B/R_u(G)$, and we take for $B^-$ the preimage in $G$ of an opposite Borel subgroup (\cite{humphreys}, \S 26.2).
Denoting by $U^-$ the unipotent part of $B^-$, we have $B \cap U^- = R_u(G)$
(as schemes), since $U^-$ contains $R_u(G)$ and the intersection
$(B/R_u(G)) \cap (U^-/R_u(G))$ is trivial. It follows that
$$
(B G_{\rm ant}) \cap U^- = B (G_{\rm ant} \cap G_{\rm aff}) \cap U^-
= B \cap U^- = R_u(G),
$$
where the first equality holds since $G_{\rm aff}$ contains $B$ and $U^-$, and
the second one since $G_{\rm ant} \cap G_{\rm aff}$ is contained in the center
of $G_{\rm aff}$ and hence in the Borel subgroup $B$.
Now the multiplication in $G$ yields a map
$$
f : B G_{\rm ant} \times U^- \to G
$$
making $B G_{\rm ant} \times U^-$ a torsor under $(B G_{\rm ant}) \cap U^- = R_u(G)$. The image of $f$ is open
in $G = G_{\rm aff} G_{\rm ant}$, since $B U^- = B B^-$ is open in $G_{\rm aff}$.
Thus the pullback of any connected \'etale cover of $G$ to the image
$B G_{\rm ant} U^-$ is connected, and hence the natural map
$$
\pi_1(B G_{\rm ant} U^-) \to \pi_1(G)
$$
is surjective. But by Lemma \ref{Gators} the map
$$
\pi_1(B G_{\rm ant} \times U^-) \to \pi_1(BG_{\rm ant} U^-)
$$
is surjective as well. This implies the surjectivity of the natural map
$$
\pi_1(B G_{\rm ant} \times U^-)^{(p')} \to \pi_1(G)^{(p')}
$$
and in turn the lemma, by using (\cite{org}, Cor. 4.9) again, and recalling
that the variety $U^-$ is an affine space.\end{dem}

\begin{lem}\label{maxtorus} The natural map
$$
\pi_1(T G_{\rm ant})^{(p')} \to \pi_1(B G_{\rm ant})^{(p')}
$$
is surjective.
\end{lem}

\begin{dem} First note that $B = U T$, and hence
$B G_{\rm ant} = U T G_{\rm ant}$. This yields a torsor
$U \times T G_{\rm ant} \to B G_{\rm ant}$
under the group scheme $U \cap T G_{\rm ant}$. The latter is commutative,
finite (since its reduced neutral component is unipotent and diagonalizable),
and unipotent, hence of order a power of $p$.
Write $G_1$ for the quotient of
$U\times TG_{\rm ant}$ by the reduced part of $U\cap G_{\rm ant}$.
The natural map
$G_1\to BG_{\rm ant}$ is finite, surjective and purely inseparable, therefore
it induces an isomorphism on fundamental groups
(\cite{sga1}, expos\'e IX, th\'eor\`eme 4.10). We may thus assume $G_1=BG_{\rm ant}$,
in which case $U\times TG_{\rm ant}\to BG_{\rm ant}$ is a finite  \'etale cover
with group $\Gamma:=U\cap TG_{\rm ant}$.
We thus obtain an exact sequence
$$
1 \to \pi_1(U \times T G_{\rm ant}) \to \pi_1(B G_{\rm ant}) \to \Gamma \to 1
$$
of fundamental groups. But $\Gamma$ is a $p$-group, and hence the natural map
$\pi_1(U \times T G_{\rm ant})^{(p')} \to \pi_1(B G_{\rm ant})^{(p')}$
is an isomorphism. We now conclude as in the proof of the previous lemma.
\end{dem}

Combining the two lemmas above we obtain the surjectivity of the natural map
$$
\pi_1(T G_{\rm ant})^{(p')} \to \pi_1(G)^{(p')},
$$
as stated in Proposition \ref{tgant}.

\begin{rema}\rm
In the course of the above proof we have in fact shown that in characteristic $p>0$ the subgroup $T G_{\rm ant}\subset G$
is a semi-abelian variety which is an extension by $T$ of an abelian variety isogenous to $A$.
If $p=0$ the same holds for the image of $T G_{\rm ant}$ in $G/R_u(G)$. A similar argument shows that (in all characteristics) the subgroup $B G_{\rm ant}$  is an extension by $B$ of an abelian variety isogenous to $A$.

One shows easily that $BG_{\rm ant}$ is a maximal connected solvable subgroup of $G$
and moreover all such subgroups are conjugate in $G$. Likewise, $T G_{\rm ant}$
(resp. the image of $T G_{\rm ant}$ in $G/R_u(G)$) is a maximal semi-abelian variety
if $p>0$ (resp. $p =0$), and all such subgroups are conjugate. Thus $B G_{\rm ant}$ is a
subgroup of $G$ that is analogous to a Borel subgroup in the linear case,
and $T G_{\rm ant}$ is analogous to a maximal torus.
\end{rema}

\section{Application to Brauer groups}\label{secberkex}

In this section we prove Theorem \ref{berkex}. We need two lemmas, in which all cohomology groups are taken with respect to the \'etale topology.

\begin{lem}\label{wedge} Let $G$ be a connected commutative algebraic group over an algebraically closed field $k$. For $n$ prime to the characteristic the cohomology ring $H^*(G, \Z/n\Z)$ with its cup-product structure is canonically isomorphic to the exterior algebra $\Lambda^*H^1(G, \Z/n\Z)$.
\end{lem}

\begin{dem}
The group $G$ is an extension of
a semi-abelian variety $G_1$ by a connected commutative unipotent group $U$. Since $U$
is an iterated extension of $\Ga$'s, we may view $G$ as an affine bundle over $G_1$.
But then the cohomology of $G$ with finite coefficients identifies with that of $G_1$.
So we may assume that $G$ is a semi-abelian variety, extension of an abelian variety $A$
by a torus $T$.
In this case the proof is a direct generalization of the argument for abelian varieties as in \cite{milneAV} Theorem 15.1. Here is a brief sketch. By formal cohomological arguments as in {\em loc. cit.} one sees that it is enough to prove the lemma with $\Q_\ell$-coefficients instead of $\Z/n\Z$-coefficients, where $\ell$ is a prime invertible in  $k$. Then one invokes an algebraic lemma (\cite{milneAV}, Lemma 15.2) based on Borel's fundamental structure theorem for Hopf algebras. It shows that it is enough to check that there exists an integer $m>0$ such that $H^i(G, \Q_\ell)=0$ for $i>m$ and moreover
$H^1(G, \Q_\ell)$ has dimension $\leq m$ over $\Q_\ell$.
We check these properties for $m=2g+r$, where $g = \dim A$ and $r = \dim T$.
In the Leray spectral sequence for $\alpha: \, G \to A$
$$
E_2^{pq}=H^{p}(A, R^q \alpha_* \Q_\ell)\Rightarrow H^{p+q}(G, \Q_\ell)
$$
the terms $E_2^{pq}$ vanish for $p>2g$ (because $A$ has  cohomological dimension $2g$) and $q>r$ (because $T$ has cohomological dimension $r$, being affine of dimension $r$). This shows the first property. For the second property, observe first that there is an isomorphism ${}_nG\cong (\Z/n\Z)^{2g+r}$ because the sequence
$$
0\to {}_nT\to {}_nG\to {}_nA\to 0
$$
is split exact by $n$-divisibility of $T$; moreover, there are isomorphisms
${}_nT\cong (\Z/n\Z)^{r}$ and ${}_nA\cong (\Z/n\Z)^{2g}$. Thus it suffices to apply the next lemma.
\end{dem}

\begin{lem}\label{sl} Under the assumptions of Lemma \ref{wedge} the map $$\Hom({}_nG, \mu_n)\to H^1(G, \mu_n)$$ sending a character $\chi:\, {}_nG\to\mu_n$ to the class of the $\mu_n$-torsor over $G$ obtained by pushing out the extension \begin{equation}\label{next} 0\to\,{}_nG\to G\to G\to 0
\end{equation} by $\chi$ is an isomorphism.
\end{lem}

\begin{dem}
 To prove injectivity, assume that a pushout of (\ref{next}) by a nontrivial character gives a split torsor over $G$. By construction it is then a split extension of $G$ by a cyclic group of order dividing $n$ which is also a quotient of $G$. This is impossible because $G$ is $n$-divisible.

 For surjectivity assume given a $\mu_n$-torsor $\pi:\,Y\to G$. Assume first $Y$ is connected. By Proposition \ref{pi1p} $a)$
there is an algebraic group structure on $Y$
making it a central extension of $G$; in particular it is commutative. (Indeed, since $\mu_n\subset Y$ is central and $G$ is commutative, the commutator map $Y\times Y\to Y,\, (y,z)\mapsto yzy^{-1}z^{-1}$ induces a map $G\times G\to\mu_n$ which must be constant as $G$ is connected.) By the same argument as for abelian varieties (\cite{mumford}, Remark p. 169 or \cite{fg}, Corollary 5.6.9) the multiplication-by-$n$ map $G\to G$ factors through $\pi$, which gives in particular a character $\chi:\, {}_nG\to\mu_n$ by looking at the fibres above the neutral element. There is thus a map of $\mu_n$-torsors from the pushforward of (\ref{next}) by $\chi$ to $Y$ which must be an isomorphism. If $Y$ is not connected, it is induced from a connected $\mu_m$-torsor for some $m|n$ which in turn comes from a character ${}_mG\to\mu_m$. The character corresponding to $Y$ is obtained by extending this character to ${}_nG$ in a way compatible with the induction of torsors.
\end{dem}

\begin{rema}\rm It follows from the above proof that the extension (\ref{next}) gives the largest abelian \'etale Galois cover of $G$ of exponent $n$. Indeed, an easy argument using the Leray spectral sequence as in the proof of Lemma \ref{wedge} shows that $H^1(G, \mu_n)$ is finite. As before, after choosing an isomorphism $\mu_n\cong \Z/n\Z$ we may identify connected $\mu_n$-torsors and \'etale $\Z/n\Z$-Galois covers. So we infer that there is a largest finite \'etale Galois cover whose Galois group is abelian of exponent $n$,
and this Galois group is isomorphic to $H^1(G, \mu_n)$. By the above proof it must be a quotient of the \'etale  cover given by (\ref{next}), but then the quotient map is an isomorphism by maximality.
\end{rema}

\noindent {\em Proof of Theorem \ref{berkex}.} The Kummer sequence in \'etale cohomology coming from the exact sequence
$$
1\to \mu_n\to \G\to \G\to 1
$$
of \'etale sheaves yields an exact sequence
$$
0\to \pic G/n\pic G \to H^2(G, \mu_n)\to {}_n\br G\to 0
$$
which we may rewrite using Lemma \ref{wedge} as
$$
0\to \pic G/n\pic G \to\Lambda^2 H^1(G, \mu_n)\otimes\mu_n^{\otimes -1}\to {}_n\br G\to 0.
$$
We may then replace $H^1(G, \mu_n)$ by $\Hom({}_nG, \mu_n)$ using Lemma \ref{sl}.\enddem

\noindent {\em Proof of Corollary \ref{corberkex}.} Setting $n=\ell^m$ in the exact sequence of Theorem \ref{berkex} and passing to the direct limit over $m$ we obtain an exact sequence
$$
0\to \pic G\otimes\Q_\ell/\Z_\ell \to \Hom(\Lambda^2 T_\ell(G), \Q_\ell/\Z_\ell)\to \br G\{\ell\}\to 0
$$
where $T_\ell(G)$ is the $\ell$-adic Tate module of $G$ and $\br G\{\ell\}$ denotes the $\ell$-primary torsion of $\br G$. As the Hom-group in the middle is $\ell$-divisible, so is $\br G\{\ell\}$, hence it is isomorphic to a finite direct power of $\Q_\ell/\Z_\ell$. To determine this power explicitly, we calculate the left and middle terms.

For the left term, we use first that $\pic G$ is isomorphic to $\pic G/U$
as in the proof of Lemma \ref{wedge}, where $U$ denotes the largest connected unipotent
subgroup of $G$. So we may assume that $G$ is an extension of $A$ by the torus $T=\G^r$.
Now recall that by virtue of a generalization of the Barsotti--Weil formula (\cite{gacc}, VII, \S 3) the extension $G$ of $A$ by $T$ is classified by a homomorphism $c:\,X^*(T)\to A^*$, where $A^*=\pic^0 A$ is the dual abelian variety and $X^*(T)$ the character group of $T$. It fits into an exact sequence
\begin{equation}\label{mag}
{\rm X}^*(T) \stackrel{c}\to \pic A \stackrel{\alpha^*}\to \pic G \to 0
\end{equation}
whose exactness follows, for instance, from (\cite{magid}, Theorem 5).

Thus $\pic G$ is the quotient of $\pic A$ by
a subgroup of $A^*$. Since $A^*$ is $\ell$-divisible, we get that $\pic G$ is an extension of
$\ns A$ by an $\ell$-divisible group, whence an isomorphism
$\pic G\otimes\Q_\ell/\Z_\ell\cong \ns A\otimes\Q_\ell/\Z_\ell$.
Moreover, we have $\ns A\cong \Z^\rho$ by definition of $\rho$.
Finally, the group $\Hom(\Lambda^2 T_\ell(G), \Q_\ell/\Z_\ell)$ is isomorphic to $ (\Q_\ell/\Z_\ell)^{(2g+r)(2g+r-1)/2}$. This follows from the fact that the $\Z_\ell$-rank of $T_\ell(G)$ is $2g+r$, as calculated at the end of the proof of Lemma \ref{wedge}.
\enddem

In conclusion, we remark that in the above proof (more precisely, in Lemma \ref{sl}) we used Proposition \ref{pi1p} in a crucial way, which in turn was based on the difficult noncommutative result of Orgogozo \cite{org}. However, for the proof of Lemma \ref{sl} it is sufficient to know that $\mu_n$-torsors over commutative algebraic groups carry a commutative group scheme structure. This fact can be proven in a much more elementary way, as we now explain.

We reduce as above to the case where $G$ is a semi-abelian variety.  A $\mu_n$-torsor $\tilde G$ over $G$ comes from a line bundle $L$ of order $n$ on $G$ together with a trivialization of $L^{\otimes n}$ (see e.g. \cite{milne}, p. 125). It thus corresponds to a class of order $n$ in $\pic G$ which, by exact sequence (\ref{mag}), comes from a class in $\pic A$. As the N\'eron--Severi group of $A$ is torsion free, the latter must be a class in $\pic^0 A$, and therefore it is represented by a line bundle whose associated $\G$-torsor $Y$ carries a group scheme structure making it a central extension of $A$ by $\G$. The $\G$-torsor $L^*$ associated with $L$ is the pullback of $Y$ to $G$, so it comes from a central extension of $G$ by $\G$. Since
$\mathcal{O}(L^*) \cong \oplus_{n = - \infty}^{\infty} H^0(G,L^{\otimes n})$, the trivialization of $L^{\otimes n}$ can be identified with a regular invertible function $\psi$ on $L^*$,  and then
$\tilde G$ is the zero subscheme of the regular function $\psi - 1$ on $G$.
When choosing the group structure on $L^*$ we may assume that $\psi$ takes the neutral element of $L^*$ to 1. But then by Rosenlicht's lemma
(\cite{rostor}, Theorem 3) $\psi$ is a character of the group $L^*$,
hence $\tilde{G} \subset L^*$ is a subgroup scheme. This defines a commutative group scheme structure on $\widetilde G$ making it a central extension of $G$ by $\mu_n$.

\end{document}